\documentclass[12pt]{amsart}
\usepackage{pb-diagram} 
\theoremstyle{plain}
\newtheorem{theorem}{Theorem}[section]
\newtheorem{prop}[theorem]{Proposition} 
\newtheorem{lemma}[theorem]{Lemma}
\newtheorem{cor}[theorem]{Corollary} 
\theoremstyle{remark}
\newtheorem{remark}[theorem]{Remark} 
\newtheorem{example}{Example}
\newcommand\lra{\longrightarrow}
\newcommand\PP[1][3]{\mathbb P^{#1}}

\newcommand\mt{\mathcal T}

\DeclareMathOperator{\tr}{tr} 
\DeclareMathOperator{\Frob}{Frob}
\begin{document}
\title[Modularity of non-rigid double octic Calabi--Yau]
{Modularity of some non--rigid double octic Calabi--Yau threefolds}
\author{S\l awomir Cynk}
\thanks{Partially supported by DFG Schwerpunktprogramm 1094
(Globale Methoden in der komplexen Geometrie) and KBN grant no. 2 P03A 013 22.}
\address{Instytut Matematyki\\Uniwersytetu Jagiello\'nskiego\\
  ul. Reymonta 4\\30--059 Krak\'ow\\Poland}
\curraddr{Institut f\"ur Mathematik\\ Universit\"at Hannover\\
Welfengarten~1\\ D--30060 Hannover\\ Germany}
\email{s.cynk@im.uj.edu.pl}
\author{Christian Meyer}
\address{Fachbereich Mathematik und Informatik\\ Johannes
  Gutenberg-Uni\-ver\-sit\"at\\
  Staudingerweg 9\\D--55099 Mainz\\Germany}
\email{cm@mathematik.uni-mainz.de}
\maketitle
\section*{Introduction}
\label{sec:intro}
The modularity conjecture for Calabi--Yau manifolds predicts that
every Calabi--Yau manifold should be modular in the sense that its
$L$--series coincides with the $L$--series of some automorphic
form(s). The case of rigid Calabi--Yau threefolds was (almost) 
solved by Dieulefait and Manoharmayum in \cite{DM,Dieulefait}.
On the other hand in the non--rigid case it is even not clear which
automorphic forms should appear. 

Examples of non--rigid modular Calabi--Yau threefolds were
constructed by Livn\'e and Yui (\cite{LivneYui}), Hulek and Verrill
(\cite{HulekVerrill, HulekVerrill2}) and Sch\"utt (\cite{Schuett}). 
In these examples modularity means a decomposition of the associated Galois
representation into two-- and four--dimensional subrepresentations with
$L$--series equal to $L(g_{4},s)$, $L(g_{2},s-1)$ or
$L(g_{2}\otimes g_{3},s)$, where $g_{k}$ is a weight $k$ cusp form.
The summand with $L$--series equal to $L(g_{2}\otimes g_{3},s)$
is explained by a double cover of a product of a K3 surface and
an elliptic curve (see \cite{LivneYui}). 

The $L$--series $L(g_{2},s-1)$ is the $L$--series of the product of
the projective line $\PP[1]$ and an elliptic curve $E$ with
$L(E,s)=L(g_{2},s)$. A two--dimensional subrepresentation with such an
$L$--series may be identified  by a map $\PP[1]\times E\lra X$ which
induces a non--zero map on the third cohomology (see \cite{HulekVerrill}).
Using an interpretation in terms of deformation theory we conjecture
that a splitting of the Galois action into two--dimensional pieces can
happen only for isolated elements of any family of Calabi--Yau threefolds. 

In this paper we will study modularity of some non--rigid double
octic Calabi--Yau threefolds, we will prove modularity of all examples listed in
table~\ref{tab1} except $X_{154}$. We will use four methods for proving
modularity, apart from the methods of Livn\'e--Yui and Hulek--Verrill
we will use two others based on giving a correspondence with a rigid Calabi--Yau
threefold or on an involution. We also observe that the splitting of the
Galois action into two--dimensional pieces holds for those Calabi--Yau
threefolds in the studied families having some additional geometric
poperties. The Calabi--Yau threefold  $X_{154}$ is also the
only one which we were not able to represent as a Kummer fibration
associated to a fiber product of elliptic fibrations
(cf. \cite{Schoen}).

\section{Modular double octics with $h^{12}=1$}

Let $D$ be an arrangement of 8 planes in $\PP$. If no six of the
planes intersect in a point and no four in a line then the double covering of
$\PP$ branched along $D$ admits a resolution of singularities $X$ which is
a smooth Calabi--Yau threefold (see \cite{Cynk}). The resolution of
singularities is performed by blowing up singularities of the branch
locus in the following order: fivefold points,
fourfold points that do not lie on a triple line,
triple lines, double lines. The Euler number of the resulting Calabi--Yau
threefold can easily be expressed in numbers of different types of
singularities. The Hodge number $h^{1,2}(X)$ (the dimension of
the deformation space) can be computed as the dimension of the
space of equisingular deformations of $D$ in $\PP$; it can also be
computed as the dimension of the equisingular ideal of $D$ (see
\cite{CvS}).  

An extensive computer search in \cite{Meyer} produced 18 double octic
Calabi--Yau threefolds with $h^{12}=1$ (in 11 one-parameter families)
for which 
\[
\tr(\Frob_{p}^{*}|H^{3}_{\text{\'et}}({X}))=a_{p}+p\cdot b_{p},
\]
for all primes $5\leq p\leq 97$, where $a_{p}$ (resp. $b_{p}$) are the coefficients of a
weight four (resp. two) cusp form.  This is a strong numerical evidence
for modularity in the sense of splitting into two two--dimensional subrepresentations.
We list all these examples in table~\ref{tab1}. We include the no. of
the arrangement (as in \cite{Meyer}), the equation, the expected modular form
of level 4 and 2 (using W. Stein's notation from \cite{modForms}) and the
Picard number $h^{11}$. Since the Calabi--Yau threefolds in the table
coming from arrangements with the same no. are birational (see
Lemma~\ref{lem:birat}) we will use in this paper
the notation $X_{n}$ for any Calabi--Yau threefold in the table constructed
from arrangement no. $n$. 

\begin{table}[htp]
  \centering \def\arraystretch{1.2}
\(
\begin{array}{|c|l|c|c|c|c|} 
\hline \text{no.} &\text{equation: } u^2=xyzt\cdot\dots&
\text{wt. 4} & \text{wt. 2}&h^{11}\\ 
\hline \hline 
4 &\textstyle  (x+y) (y+z) (x-y-z-t) & 32k4A1 &
32A1&61\\
&(x+y-z-t)&&&\\ 
\hline 4 &  (x+y) (y+z) (x+2y+2z-t)& 32k4A1 &32A1&61\\ 
& (x+y+2z-t)&&&\\ 
\hline 
4 & 2 (x+y) (y+z) (2x+y+z-2t) & 32k4A1 & 32A1&61\\ 
&(2x+2y+z-2t)&&&\\ 
\hline 
8 &(x+y) (y+z) (-z+t) (3x-y-z+t)& 24k4A1 & 24A1&61\\ 
\hline 
13 & (x+y) (y+z) (x-z-t) (x-z-2t)& 32k4A1 & 32A1&61\\ 
\hline 
13 & (x+y) (y+z) (x-z-t) (x-z+t)& 32k4A1 & 32A1&61\\ 
\hline 
13 &  (x+y) (y+z) (x-z-t) (2x-2z-t)& 32k4A1 & 32A1&61\\ 
\hline 
21 &  (x+y) (y+z) (2x+y-t) (2x-z-2t)& 32k4B1 & 32A1& 53\\ 
\hline
53 &  (x+y) (z+t) (x-y-z-t)& 32k4B1 & 32A1& 53\\
& (x+y-z+t)&&&\\ 
\hline 
154 & (x+y+z) (x+y+z-t) &
8k4A1 & 72A1& 41 \\
& (-2x+y-3z+3t)(2x+3z-2t)&&& \\
 \hline 
244 & (x+y+z+t)(x+y-z-t)(y-z+t)&12k4A1 & 48A1& 39 \\ 
& (x-z+t)&&&\\ 
\hline 
249 & (x+y+z) (x+z+t) (2x+3y-z+2t)&24k4A1 & 24A1& 37 \\ 
&(y-z+2t)&&&\\ 
\hline 
249 & (x+y+z) (x+z+t) (2x-y+3z+2t)&24k4A1 & 24A1& 37 \\ 
& (-3y+3z+2t)&&&\\ 
\hline 
267 & (x+y-2z)(x-y-z+t)(2y-z+t)&96k4B1 & 96B1& 37 \\ 
&(x+y+z+t)&&&\\ 
\hline
267 & (x+y+z) (x+2y-z+t)& 96k4B1 & 96B1& 37 \\ 
&
(-y+2z-2t)(2x+2y-z+2t)&&&\\ 
\hline 
267 & (2x+2y-z) (2x+y-2z+2t)
&96k4B1 & 96B1& 37 \\ 
&(y+z-t)(x+y-2z+t)&&&\\ 
\hline 
274 & (x+y+z)
(-x-z+t)(x+2y-z+t) &96k4E1 & 96B1& 37 \\ 
& (x+y-z+2t)&&&\\ 
\hline
275 & (x+y+z)(2x-2z-t)(8y+4z+t)&96k4B1 & 96B1& 37 \\
&(2x+4y+t)&&&\\ \hline
\end{array} \)\\[2mm]
\caption{}
\label{tab1}
\end{table}

The Picard groups of all listed Calabi--Yau threefolds are generated
by divisors defined over $\mathbb Q$, so Frobenius acts on
$H^{2}_{\text{\'et}}$ by multiplication with $p$. In fact, in all the examples
except $X_{244}$, the skew-symmetric part of the
Picard group is zero, whereas for $X_{244}$ it is generated by a
divisor coming from the contact plane $x+y-z+t=0$.

\section{Double quartic elliptic fibrations}
\label{sec:dq}

In this section we will shortly review some information about
rational elliptic fibrations that can be realized as a resolution of a
double covering of $\PP[2]$ branched along a sum of four lines. The
structure of the elliptic fibration is determined by the choice of a point
in $\PP[2]$. Some of these surfaces where described in \cite{CynkMeyer2};
we will omit here all the details explained in that paper.

The double covering is rational exactly when the lines do not intersect in
one point. We can have the following combinations of singular fibers
(the Picard number $\rho(S_{w})$ of a generic fiber can be computed from
the Zariski lemma):
\[  \begin{array}{l|l|c}
  &\text{singular fibers}&\rho(S_{w})\\\hline
   S_{1} &D_{4}^{*},D_{4}^{*}&1\\
   S_{2} &I_{2}, I_{2}, D_{6}^{*}&1\\
   S_{3} &I_{2}, I_{2},  I_{4}, I_{4}&1\\
   S_{4} &I_{2}, I_{2}, I_{2}, D_{4}^{*}&2\\
   S_{5} &I_{2}, I_{2}, I_{2}, I_{2}, I_{4}&2\\
   S_{6} &I_{2},I_{2},I_{2},I_{2},I_{2},I_{2}&3
  \end{array}
\]

A double covering of $S_{1}$ branched along the two singular fibers is
birational to a product of $\PP[1]$ and an elliptic curve $E$, and all
smooth fibers are isomorphic to $E$. This elliptic fibration
depends on the $j$--invariant of~$E$. 

The surfaces $S_{2}$ and $S_{3}$ are extremal, i.e. they have
$\rho(S_{w})=1$. Consequently they are uniquely defined as fiber
spaces. Moreover the parameters corresponding to the singular fibers of
$S_{3}$ form a harmonic quadruple (i.e. their cross--ratio equals $-1$);
they can be chosen as
\[
\begin{array}{rrrr} 
-1&0&1&\infty\\ \hline
I_{2}&I_{4}&I_{2}&I_{4}
\end{array}
\]

Denote by $S_{3}'$ the pullback of $S_{3}$ via the involution $t\mapsto
\frac{t-1}{t+1}$ of $\PP[1]$, so $S_{3}'$ has the following singular
fibers:
\[
\begin{array}{rrrr} 
-1&0&1&\infty\\ \hline
I_{4}&I_{2}&I_{4}&I_{2}
\end{array}
\]
Thus $S_{3}$ and $S_{3}'$ have singular
fibers at the same points but of different types. There exists an isogeny
$\Psi:S_{3}\mapsto S_{3}'$ which is a degree 2 unbranched covering
on a smooth fiber. 

Fibration $S_{4}$ is not extremal, so we can chose arbitrary
coordinates of singular fibers. The configuration of lines is not
uniquely determined by the coordinates of singular fibers. In fact
there are exactly two types: one with a triple point and one with a
``vertical line''. 

The Picard number of the generic fiber of Fibration $S_{5}$ equals  two,
so we can not choose arbitrary coordinates of singular fibers. In fact
there is an involution of $\PP[1]$ which preserves the fiber $I_{4}$
and exchanges two pairs of $I_{2}$'s. The configuration of lines is
uniquely determined. 

Fibration $S_{6}$ is the most complicated one. In this case the
configuration of lines is not uniquely determined. There can be
several choices coming from automorphisms of $\PP[1]$ preserving the
singular fibers.

\section{Kummer fibrations}
\label{sec:kf}

All examples in table~\ref{tab1} except $X_{154}$ can be realized as a Kummer
fibration associated to a fiber product of elliptic fibrations (cf. \cite{Schoen}).
Contrary to Schoen we do not require that the involution on the fiber
product lifts to a resolution, so the resulting Calabi--Yau threefold
is not necessarily a blow--up of the Kummer fibration.

To see the fibration we  reorder the planes such that the first four
and the last four intersect in a point. Then after change of
coordinates in $\PP$ we may assume that these points of intersection
are $(0,0,0,1)$ and $(1,0,0,0)$, or equivalently that the double octic
is given in weighted projective space $\PP[](1,1,1,1,4)$ by the equation
\[
w^{2}=f_{1}(x,y,z)\cdot\ldots\cdot f_{4}(x,y,z)f_{5}(y,z,t)\cdot\ldots\cdot
f_{8}(y,z,t).
\]
Consequently the double octic is birational to the
quotient of the fiber product of elliptic fibrations
\[
u^{2}=f_{1}(x,y,z)\cdot\ldots\cdot f_{4}(x,y,z)
\] and
\[
v^{2}=f_{5}(y,z,t)\cdot\ldots\cdot f_{8}(y,z,t)
\]
by the involution
\[
(x,y,z,t,u,v)\mapsto(x,y,z,t,-u,-v).
\]

In the following table we list descriptions of Calabi--Yau threefolds from
table~\ref{tab1} as Kummer fibrations. For each Kummer fibration we
give coordinates and types of singular fibers. In some cases we
were able to find two different representations as a Kummer fibration.
\[\def\arraystretch{1.2}
\begin{array}{ll}
\rule{15mm}{0cm}&\rule{9cm}{0cm}\\ 
X_{4}&
\begin{array}{rrrrr} 
0&1&2&3&\infty\\ \hline
I_{2}&I_{2}&I_{2}&I_{2}&I_{4}\\
I_{0}&I_{2}&I_{2}&I_{0}&D_{6}^{*}\\
\end{array}\hfill\begin{array}{rrrr} 
-1&0&1&\infty\\\hline
I_{4}&I_{2}&I_{4}&I_{2}\\ 
D_{6}^{*}&I_{2}&I_{0}&I_{2}
\end{array}
\\
 \hline
X_{8}&
\begin{array}{rrrr} 0&1&4&\infty\\ \hline
D_{4}^{*}&I_{2}&I_{2}&I_{2}\\
I_{2}&I_{2}&I_{0}&D_{6}^{*}
\end{array}\\
 \hline
X_{13}&
\begin{array}{rrr} 0&1&\infty\\ \hline
D_{4}^{*}&D_{4}^{*}&I_{0}\\
I_{2}&I_{2}&D_{6}^{*}
\end{array}\hfill
\begin{array}{rrrr} 
-1&0&1&\infty\\ \hline
I_{2}&I_{4}&I_{2}&I_{4}\\ 
I_{0}&D_{4}^{*}&I_{0}&D_{4}^{*}
\end{array}
\\ \hline
X_{21}&
\begin{array}{rrrr} 
-1&0&1&\infty\\ \hline
I_{2}&I_{4}&I_{2}&I_{4}\\
D_{6}^{*}&I_{0}&I_{2}&I_{2}
\end{array}
\hfill
\begin{array}{rrrr} 
-1&0&1&\infty\\\hline
I_{2}&I_{2}&D_{4}^{*}&I_{2}\\ 
D_{4}^{*}&I_{2}&I_{2}&I_{2}
\end{array}
\\
 \hline
X_{53}&
\begin{array}{rrrr} 
-1&0&1&\infty\\ \hline
I_{2}&D_{6}^{*}&I_{2}&I_{0}\\
I_{2}&I_{0}&I_{2}&D_{6}^{*}
\end{array}
\hfill
\begin{array}{rrrr} 
-1&0&1&\infty\\ \hline
I_{2}&I_{2}&I_{2}&D_{4}^{*}
\\I_{2}& D_{4}^{*}&I_{2}&I_{2}
\end{array}
\\ \hline
X_{244}&
\begin{array}{rrrrr} 
-1&0&1&2&\infty\\ \hline
I_{0}&I_{2}&I_{4}&I_{2}&I_{4}\\
I_{4}&I_{2}&I_{4}&I_{0}&I_{2}
\end{array}
\hfill
\begin{array}{rrrrrr} 
-1&0&\frac13&1&3&\infty\\ \hline
I_{2}&I_{2}&I_{2}&I_{4}&I_{0}&I_{2}
\\I_{2}&I_{2}&I_{0}& I_{4}&I_{2}&I_{2}
\end{array}
\\ \hline
X_{249}&
\begin{array}{rrrrrr} 
-1&0&\frac13&1&3&\infty\\ \hline
I_{0}&I_{2}&I_{2}&I_{4}&I_{2}&I_{2}
\\I_{2}&I_{4}&I_{0}& I_{2}&I_{0}&I_{4}
\end{array}\\
\hline
X_{267}&
\begin{array}{rrrrrr} 
-1&0&\frac12&1&2&\infty\\ \hline
I_{2}&I_{2}&I_{2}&I_{2}&I_{2}&I_{2}
\\I_{2}&I_{2}&I_{2}& I_{2}&I_{2}&I_{2}
\end{array}\\
 \hline
X_{274}&
\begin{array}{rrrrrr} 
-1&0&\frac12&1&2&\infty\\ \hline
I_{2}&I_{4}&I_{2}&I_{2}&I_{0}&I_{2}
\\I_{4}&I_{2}&I_{2}& I_{0}&I_{2}&I_{2}
\end{array}\\
 \hline
X_{275}&
\begin{array}{rrrrrr} 
-1&0&\frac12&1&2&\infty\\ \hline
I_{2}&I_{2}&I_{2}&I_{2}&I_{2}&I_{2}
\\I_{2}&I_{2}&I_{2}& I_{2}&I_{2}&I_{2}
\end{array}
\end{array}
\]

\begin{lemma}\label{lem:birat}
The Calabi--Yau threefolds in table \ref{tab1} defined by arrangements of the
same type are birational. The Calabi--Yau threefolds $X_{21}$ and $X_{53}$ are
birational; and the Calabi--Yau threefolds $X_{267}$ and $X_{275}$ are birational. 
There exists a correspondence between the Calabi Yau--threefolds $X_{8}$
and~$X_{249}$. 
\end{lemma}
\begin{proof}
From the explicit description of the fiber products in local coordinates it
easily follows that the Calabi--Yau threefolds defined by arrangements of
the same type with different parameters are in fact projectively equivalent.

Arrangement no. 21 is projectively equivalent to 
\[
x(x-z)(x+z)(x+y)y(t+z)(t-z)(t+y)=0.
\]
Substituting the birational involution of $\PP$ given by
\[
(x,y,z,t)\mapsto(yz,xz,xy,tx)
\]
we obtain
\[
(xzy^{2})^{2}x(x-z)(x+z)(x+y)z(t+y)(t-y)(t+z)=0,
\]
and since arrangement no. 53 is projectively equivalent to 
\[
x(x-z)(x+z)(x+y)z(t+y)(t-y)(t+z) = 0,
\]
we conclude that the resulting Calabi--Yau threefolds are birational. 

To prove that $X_{267}$ and $X_{275}$ are birational, observe that the
corresponding arrangements are projectively equivalent to

\begin{eqnarray*}
  &\text{Arr. no. 267:}\quad&x(x-z)(2x-2z+y)(2x-z-y)\times\\
  &&\times t(t+z-y)(2y-z-2t)(2z-y+2t)=0\\
  &\text{Arr. no. 275:}\quad&x(x-z)(2x-2z+y)(2x-z-y)\times\\
  &&\times t(2t-y)(2t-z)(3t-y-z)=0.
\end{eqnarray*}
Simple computations show that the cross ratios of the quadruples
\[
\begin{array}{llll}
0, & \quad y-1, & \quad y-\tfrac12, & \quad \tfrac12y-z\\
0, & \quad \tfrac 12y, & \quad \tfrac 12z, & \quad \tfrac 13y+\tfrac 13z
\end{array}
\]
are equal so there is a birational transformation in $y,z,t$ that
maps one of them to the other.

To see the correspondence between the Calabi--Yau threefolds $X_{8}$
and $X_{249}$, first pull back arrangement no. 8 by the
map $t\mapsto (\frac{t+1}{t-1})^{2}$, obtaining 
\[
\begin{array}{rrrrrr} 
-1&0&\frac13&1&3&\infty\\ \hline
I_{0}&I_{2}&I_{2}&I_{4}&I_{2}&I_{2}
\\I_{4}&I_{2}&I_{0}& I_{4}&I_{0}&I_{2}
\end{array}\]
Now it is enough to compose this map with the isogeny of the elliptic
fibration with fibers $I_{4},I_{4},I_{2},I_{2}$ that exchanges $I_2$ fibers
with $I_4$ fibers (see \cite{CynkMeyer2}).
\end{proof}

\begin{remark}
Arrangements no. 267 and 275 are not projectively equivalent, they
come from different twisted self--fiber products of the same
elliptic fibration. The self--fiber product (without twist) of this
elliptic fibration gives a non--birational Calabi--Yau threefold
with $h^{12}=2$ (see example~\ref{ex:ar269}). 
\end{remark}

\section{Ruled surface over elliptic curves}
\label{sec:ell}

In this section we will use elliptic ruled surfaces to prove
modularity of four Calabi--Yau threefolds from table \ref{tab1}.
\begin{prop}
The Calabi--Yau threefolds $X_{4}$, $X_{8}$, $X_{244}$ and $X_{249}$ are
modular, with modular forms as listed in table \ref{tab1}.
\end{prop}

Consider a Calabi--Yau threefold $X$ such that an $L$--series of the
form $L(g_{2},s-1)$ (where $g_{2}$ is a weight two modular form corresponding
to an elliptic curve $E$) appears in the Galois representation. Then
by the Tate Conjecture we can expect that there is a correspondence
between $X$ and the product $E\times \PP[1]$ which induces the
isomorphism of representations. 

Hulek and Verrill proved in \cite{HulekVerrill} that when a smooth
ruled surface over an elliptic curve $S\lra E$ 
is contained in a Calabi--Yau threefold $X$ then the map on third
cohomology $H^{3}(X)\lra H^{3}(S)$ is surjective. The map can be
represented by a direct sum of $H^{1}(\mt _{X})\lra H^{1}(\mathcal
N_{S|X})$ and its complex conjugate. The map $H^{1}(\mt _{X})\lra
H^{1}(\mathcal N_{S|X})$ associates to a deformation of $X$ the
obstruction to lift it to a deformation of $E$ (inside $X$). Therefore
if this map is non-zero then $E$ deforms inside $X$ only over a
codimension one submanifold of the Kuranishi space of $X$.

Now, if we have ruled surfaces $E_{1},\dots,E_{r}$, with
$r=h^{21}(X)$, such that the map
\begin{equation}
  \label{eq:thirdcoho}
H^{3}(X)\lra \bigoplus_{i} H^{3}(E_{i})  
\end{equation}
is surjective then the
obstructions are independent and the surfaces do not deform
simultaneously over any subvariety of the Kuranishi space of $X$ of
positive dimension. It is an explanation why in a family there were
always only finitely many examples were one was able to prove
modularity in that way.

If we have several ruled surface over elliptic curves, it is usually
difficult to determine whether the map \eqref{eq:thirdcoho}
is surjective. In case we know the Kuranishi space of
$X$ we can try to invert the above argument. For each elliptic
fibration we consider the hypersurface $V_{i}$ of the Kuranishi space over
which $E_{i}$ deforms, knowing that the kernel of \eqref{eq:thirdcoho} is the
tangent to the intersection of the $V_{i}$'s plus its complex
conjugate (see example at the end of this section).
 
To use this method in our examples we need to find elliptic fibrations inside
the double octics. If a plane $S$ in $\PP[3]$ contains two double lines
and the other four arrangement planes intersect at a point in $S$,
then the pullback of $S$ to the double covering is an elliptic
fibration. On the Kummer fibration these planes are recognized as
corresponding to the product of fibers $I_{0}$ and $I_{4}$.

We were able to find such a plane only for two arrangements:

{\bf Arrangement no. 4:}
the plane $S$ has equation $x-z=0$ resp. $y+2z-t=0$ resp. $2x+y-2t=0$
(for the three arrangements in the table).

{\bf Arrangement no. 244:}
the plane $S$ has equation $x+y+z-t=0$.

\medskip
To prove modularity of $X_{8}$ and $X_{249}$ we will study an
auxiliary Calabi--Yau threefold $X_{269}$ with $h^{12}=2$. Modularity
of this Calabi--Yau threefold follows from existence of some elliptic ruled
surfaces and their behavior under deformations.

\begin{example}\label{ex:ar269}
Consider the double octic Calabi--Yau threefold $X_{269}$ defined by the
following arrangement of eight planes (arrangement no. 269 in \cite{Meyer}):
\[
xyzt(x+y+z)(x+2y-z+t)(y+z-t)(x+y-2z+t)=0
\]
It has $h^{2,1}(X_{269})=2$.
Substituting $y=y-z,z=z+t$ we can represent this Calabi--Yau
threefold as the following Kummer fibration:
\[
\begin{array}{rrrrrr} 
-1&0&\frac13&1&3&\infty\\ \hline
I_{2}&I_{2}&I_{2}&I_{2}&I_{2}&I_{2}
\\I_{2}&I_{2}&I_{2}& I_{2}&I_{2}&I_{2}
\end{array}\]
On the other hand substituting $x=x+2z-4y, z=x-2y$ we can also obtain
the following Kummer fibration:
\[
\begin{array}{rrrrrr} 
-1&0&\frac13&1&3&\infty\\ \hline
I_{0}&I_{2}&I_{2}&I_{4}&I_{2}&I_{2}
\\I_{4}&I_{2}&I_{0}& I_{4}&I_{0}&I_{2}
\end{array}
\]

Hence using the isogeny between $S_{3}$ and $S_{3}'$ from section~\ref{sec:dq} we can find
correspondences between this Calabi--Yau threefold and the Calabi--Yau
threefolds $X_{8}$ and $X_{249}$. 

Observe that the planes $z=x+2y$ and $y=2z-t$ contain two double
lines and a fourfold point, so they give two ruled surfaces $E_{1},
E_{2}$ over an 
elliptic curve with conductor 24. 

The Kuranishi space of the Calabi--Yau threefold $X_{269}$ may be
para\-met\-rized by the equation
\begin{align*}
& xyzt(x+y+z)(Bx+Cy-Az+At) \,\times\\
&\qquad \times (y+z-t)(Bx+By+(-A+B-C)z+At)=0.      
\end{align*}
By \cite{HulekVerrill2} both elliptic fibrations give non-zero maps  
\[H^{3}(X)\lra H^{3}(E_{i})\]
so they deform over curves in $\PP[2]$.
One easily checks that they deform over the lines given by
\begin{align*}
A+B-C&=0,\\C&=2B,
\end{align*}
which intersect only at the point $(1,1,2)$ corresponding to the
equation we started with. Consequently the obstructions are
independent and the map
\[H^{3}(X)\lra H^{3}(E_{1})\oplus H^{3}(E_{2})\]
is surjective, giving a splitting of the representation on $H^{3}$
into two--dimensional pieces. Counting points over $\mathbb F_{p}$
for $p\le97$ one checks that $X$ is modular and that the coefficients of
the $L$--series are given by $b_{p}+2pc_{p}$, where $b_{p}$ resp. $c_{p}$
are the coefficients of the unique cusp form of level 24 and weight 4 resp. 2.
\end{example}

There is a degree two correspondence between the above Calabi--Yau
threefold and $X_{249}$, hence also $X_{8}$. These correspondences prove
the modularity of $X_{8}$ and $X_{249}$.

\section{Correspondences with rigid double octics}
\label{sec:rigid}

In this section we will use correspondences between rigid and
non--rigid Calabi--Yau threefolds to prove modularity of the latter.
\begin{prop}
\label{prop_4_21_53_244}
  The Calabi--Yau threefolds $X_{4}$, $X_{21}$, $X_{53}$ and $X_{244}$ are
  modular, with modular forms as listed in table \ref{tab1}.
\end{prop}

In \cite{CynkMeyer} we checked the modularity and computed modular forms
of some rigid double octic Calabi--Yau threefolds. Now we will use
correspondences between some rigid and non--rigid Calabi--Yau
threefolds to show the modularity of the latter.

We first recall the considered rigid examples. As before we
will use the equations and numbers of arrangements from
\cite{Meyer} (in brackets we give the numbers from \cite{CynkMeyer}).

\textbf{Arrangement no. 3} (old no. 6) is given by the equation
\[
xyzt(x+y)(y+z)(z+t)(t+x)=0.
\]
The corresponding fiber product of elliptic fibrations has singular fibers
\[\begin{array}[c]{cccc}
I_{4}&I_{4}&I_{2}&I_{2}\\
D_{6}^{*}&I_{2}&I_{2}&I_{0}
\end{array}
\]

\textbf{Arrangement no. 19} (old no. 23) is given by the equation
\[
xyzt(x+y)(y+z)(x-z-t)(x+y+z-t)=0.
\]
The corresponding fiber product of elliptic fibrations has singular fibers
\[
\begin{array}[c]{cccc}
I_{2}&I_{2}&I_{4}&I_{4}\\
I_{0}&D_{6}^{*}&I_{2}&I_{2}
\end{array}
\]

\textbf{Arrangement no. 239} (old no. 86$^{a}$) is given by the equation
\[
xyzt(x+y+z)(x+y+t)(x+z+t)(y+z+t)=0.
\]
The corresponding fiber product of elliptic fibrations has singular fibers
\[
\begin{array}[c]{ccccc}
I_{2}&I_{2}&I_{4}&I_{4}&I_{0}\\
I_{0}&I_{4}&I_{2}&I_{4}&I_{2}
\end{array}
\]

\begin{lemma}
  There are correspondences between the Calabi--Yau threefolds given by the
  following arrangements:
  \begin{enumerate}
  \item No. 4 and no. 19,
  \item No. 21 and no. 3,
  \item No. 53 and no. 3,
  \item No. 244 and no. 239.
  \end{enumerate}
\end{lemma}

\begin{proof}
  All the correspondences are in fact defined on the level of the fiber
  products of elliptic fibrations. They are given by applying the
  isogeny of the elliptic fibration with fibers $I_{2},I_{2},I_{4},I_{4}$
  that exchanges the fibers $I_{2}$ and $I_{4}$. 
\end{proof}

Assume that we have a generically finite correspondence between two
Calabi--Yau threefolds $X$ and $Y$. Then this correspondence induces an
isomorphism between $H^{3,0}(X)$ and $H^{3,0}(Y)$ coming from a pullback
of the canonical form. If $Y$ is rigid then taking this isomorphism plus its
complex conjugate we obtain a splitting of the Galois representation on
$H^{3}(X)$ into a two--dimensional representation isomorphic to
$H^{3}(Y)$ and its complement. Using the correspondences from the
above lemma and counting points in $\mathbb F_{p}$ for $p\le97$ we
obtain proposition \ref{prop_4_21_53_244}.

\section{Kummer construction}
\label{sec:mr}

In this section we will use the Kummer construction studied by Livn\'e and
Yui (\cite{LivneYui}).

\begin{prop}
The Calabi--Yau threefold $X_{13}$ is modular, with modular forms as listed
in table \ref{tab1}.
\end{prop}

We will consider a two--dimensional family of double octic
Calabi--Yau threefolds which are the quotient by an involution of a product
of a K3 surface studied in \cite{AOP} and an elliptic curve.
Take the elliptic curve
\[
E_{\mu}=\{(x,t,u)\in\PP[](1,1,2):u^{2}=(x-t)(x^{2}-\mu t^{2})t \}
\]
and the K3 surface
\[
S_{\lambda}=\{(y,z,t,v)\in\PP[](1,1,1,3): v^{2}=yzt(y+t)(z+t)(y+\lambda z)\}.
\]

On the product $Y_{\lambda,\mu}:=E_{\mu}\times S_{\lambda}$ we have a natural involution
\[
((x,t,u),(y,z,t,v))\longmapsto ((x,t,-u),(y,z,t,-v)).
\]

The quotient $X_{\lambda,\mu}$ of $Y_{\lambda,\mu}$ by this involution has a
Calabi--Yau nonsingular model. To show this observe that $Y_{\lambda,\mu}$
is birational to the double covering of $\PP$ branched along the octic
$D_{\lambda,\mu}$ given by the equation
\[
(x-t)(x^{2}-\mu t^{2})yz(y+t)(z+t)(y+\lambda z)=0.
\]

The birational map can be given by in appropriate affine
coordinates ($t=1$) by
\[
(x,1,u),(y,z,1,v)\longmapsto (x,y,z,uv).
\]

The octic itself is defined over $\mathbb Q$. Over $\mathbb Q[\sqrt\mu]$
it splits into a sum of eight planes (for general $\mu$, two of them are not defined
over $\mathbb Q$). Using \cite{Cynk} we conclude that $X_{\lambda,\mu}$ has a
nonsingular model $\tilde X_{\lambda,\mu}$ which is a Calabi--Yau threefold.

For general values of $\lambda$ and $\mu$, the arrangement $D_{\lambda,\mu}$
is arrangement no. 52 in \cite{Meyer}, so $\tilde X_{\lambda,\mu}$ has the
invariants $h^{11}(\tilde X_{\lambda,\mu})=56$ and $h^{12}(\tilde X_{\lambda,\mu})=2$.

For  $\lambda\not=0,-1$, the rank of the symmetric part of the Picard
group of the K3 surface $S_{\lambda}$ is 19; denote by
$H^{2}_{skew}(S)$ the three-dimensional skew-symmetric part.
Thus there is a Shioda--Inose structure on $S_{\lambda}$, namely there exists
an involution on $S_{\lambda}$ such that the quotient of $S_{\lambda}$ by that
involution is a Kummer surface.

In \cite {AOP} it is proved that the surface $S_{\lambda}$, with $\lambda\in
\mathbb Q\setminus \{0,-1\}$, is modular exactly when
$\lambda\in\{1,8,1/8,-4,-1/4,-64,-1/64\}$, and the modular form for $S_{\lambda}$
is computed.  We have the following diagram of rational maps
\[
\begin{diagram} \node{Y_{\lambda,\mu}}
\arrow{se}\arrow[2]{e,..}\node[2]{\tilde X_{\lambda,\mu}}\arrow{sw}\\
\node[2]{X_{\lambda,\mu}}
\end{diagram}
\]
The rational map $Y_{\lambda,\mu}\longrightarrow\tilde X_{\lambda,\mu}$ can
be resolved by blowing up at points and lines so it induces a well defined map in
cohomologies $H^{3}(X_{\lambda,\mu})\lra H^{3}(Y_{\lambda,\mu})$. The image of the
map is invariant under the involution on $Y_{\lambda,\mu}$, so in fact we
obtain a map $H^{3}(X_{\lambda,\mu})\lra H^{1}(E_{\mu})\otimes
H^{2}_{skew}(S_{\lambda})$. From the description of deformations of double coverings of
smooth algebraic varieties (\cite{CvS}) it follows that this map is surjective,
moreover both vector spaces have dimension 6, so it is an
isomorphism. We obtain

\begin{prop}
$H^{3}(\tilde X_{\lambda\mu})\cong H^{1}(E_{\mu})\otimes H^{2}_{skew}(S_{\lambda})$.
\end{prop}
\begin{cor} The Calabi--Yau threefold $\tilde X_{\lambda,\mu}$ is modular for
$\lambda\in\{1,8,1/8,-4,-1/4,-64,-1/64\}$ and $\mu\in\mathbb Q\setminus\{0,1\}$.
\end{cor}

For the seven values of $\lambda $ the $L$--series of $S_{\lambda}$
corresponds to a cusp form for 
$S_{3}(\Gamma_{1}(8))$, $S_{3}(\Gamma_{1}(16))$, $S_{3}(\Gamma_{1}(12))$,
$S_{3}(\Gamma_{1}(7))$ (for $\lambda$ and $1/ \lambda$ the $L$--series differ
only by a twist). They are the only $\eta$--product weight 3 modular
forms. The modular form of the surface $S_{\lambda}$ corresponds to the
symmetric power of the modular form associated to the elliptic curve
$E_{\frac 1{\lambda +1}}$ (see \cite{AOP}). For the seven special values of
$\lambda$ the elliptic curve $E_{\frac1{\lambda+1}}$ has complex multiplication.
Denoting by $a_{p}$ resp. $b_{p}$ the Fourier coefficients of the level 2
(resp. level 3) modular forms we get
\[
b_{p}=  \begin{cases}
    a_{p}^{2}-2p,\qquad&\left(\frac{-(\lambda+1)}p\right)=1\\[3mm]
    0,&\left(\frac{-(\lambda+1)}p\right)=-1.
  \end{cases}
\]
The Fourier coefficient of the $L$--series of $S_{\lambda}$ equals $\left(\tfrac{-(\lambda+1)}p\right)(b_{p}+p)$.

The third symmetric power of a weight 2 form yields also a weight 4
modular form with Fourier coefficients
\[
c_{p}=a_{p}^{3}-3pa_{p},
\]
so we obtain
\[
a_{p}b_{p}=c_{p}+pa_{p}.
\]

Consequently we get much better modularity properties for the threefolds
$X_{\lambda}:=\tilde X_{\lambda,\frac{1}{\lambda+1}}$.

\begin{prop} The $L$--series of the Calabi--Yau threefold $X_{\lambda}$ has
Fourier coefficients equal to
\[
c_{p}+2pa_{p}.
\]
\end{prop}

In the table we collect the data for the four Calabi--Yau threefolds
the $L$--series of which do not only differ by a twist:

\[
\begin{array}[t]{l||c|c|c|c} &\lambda=1&\lambda=8&\lambda=-4&\lambda=-64\\ \hline\hline
\text{wt 2 form}&256k2D&32k2A&144k2B&49k2A\\ \hline \text{wt 3
form}&8k3A[1,1]&16k3A[1,0]&12k3A[0,1]&7k3A[3]\\ \hline \text{wt 4
form}&256k4H&32k4A&144k4A&49k4D\\ \hline
b_{p}=a_{p}^{2}-2p&p\equiv1,3(8)&p\equiv3(4)&p\equiv1(3)&p\equiv1,2,4(7)\\ \hline
b_{p}=0&p\equiv5,7(8)&p\equiv1(4)&p\equiv2(3)&p\equiv3,5,6(7)\\ \hline
\parbox{2cm}{$\eta$--products\\(wt 3)}&
\parbox{22.5mm}{$\eta^{2}(z)\eta^{2}(2z)\\[1mm]\eta^{2}(4z)\eta^{2}(8z)$}&
\eta^{6}(4z)&\eta^{3}(2z)\eta^{3}(6z)&\eta^{3}(z)\eta^{3}(7z)\\ \hline
\parbox{2cm}{$\eta$--products\\(wt 2)}&-&\eta^{2}(8z)\eta^{2}(4z)&
\frac{\eta^{12}(12z)}{\eta^{4}(24z)\eta^{4}(6z)}&-
\end{array}
\]

\subsection{Singular K3}

From the above considerations we excluded the case of $\lambda=-1$. There
are two reasons for this. First, in this case all divisors on the K3 surface are
symmetric and consequently $h^{12}(\tilde X_{-1,\mu})=1$ (this is
arr. no 13). Second, $\frac 1{\lambda+1}$ makes no sense. We can however
take in that case also the curve $E_{1/9}$, as the modular forms
appearing in $S_{-1}$ and $S_{8}$ are the same. Hence for the Calabi--Yau
threefold $\tilde X_{-1,1/9}$ the modular form has coefficients
$c_{p}+pa_{p}$, where $c_{p}$ resp. $a_{p}$ are coefficients of a weight
4 resp. 2 level 32 newform.

In the above considerations we can replace the elliptic curve
$E_{\frac1{\lambda+1}}$ by another elliptic curve with the same modular
form, or replace both $E_{\mu}$ and $S_{\lambda}$ by some twist.

Now fix $\lambda\in\{1,8,\frac 18,-4,-\frac 14,-64,-\frac 1{64}\}$. 
Using \cite{AOP}
we can compute the characteristic polynomial of
Frobenius on $H^3$ for the Calabi--Yau threefold $\tilde X_{\lambda,\mu}$
for any rational $\mu\not=0,-1$.  Denoting by $\alpha_{p},\bar\alpha_{p}$ resp.
$\beta_{p},\bar\beta_{p}$ the eigenvalues of Frobenius on
$H^{1}(E_{\lambda,p})$ resp. $H^{1}(E_{\mu,p})$ we find
that the characteristic polynomial of Frobenius acting on
$H^{3}(\tilde X_{\lambda,\mu})$ is (up to sign)
\[
(T-p\beta_{p})(T-p\bar\beta_{p})\cdot(T-\alpha^{2}_{p}\beta_{p})(T-\alpha^{2}_{p}\bar\beta_{p})
(T-\bar\alpha^{2}_{p}\beta_{p})(T-\bar\alpha^{2}_{p}\bar\beta_{p}).
\]
This polynomial splits over $\mathbb Z$ into the characteristic
polynomial of the Frobenius action on $H^{2}((\PP[1]\times E_{\mu})_{p})$ and
the degree 4 polynomial $(T-\alpha^{2}_{p}\beta_{p})(T-\alpha^{2}_{p}\bar\beta_{p})
(T-\bar\alpha^{2}_{p}\beta_{p})(T-\bar\alpha^{2}_{p}\bar\beta_{p})$.
In the construction, this splitting comes from the cartesian product
of $E_{\mu}$ and a transcendental cycle on the K3 surface $S_{\mu}$;
it should have a better geometric interpretation via the Shioda--Inose structure.

If the elliptic curves $E_{\lambda}$ and $E_{\mu}$ are
non--isogenous,  the degree 4 polynomial does not divide by the characteristic
polynomial of $\PP[1]\times E$, for any elliptic curve $E$. To see this,
denote the eigenvalues of Frobenius on $H^{1}(E)$ by $\gamma _{p},\bar
\gamma_{p}$ and assume that $p\gamma_{p}=\bar\beta_{p}\alpha_{p}^{2}$.
Multiplying by $\beta_{p}$ and dividing by $p=|\beta_{p}|^{2}$ we get
$\beta_{p}\gamma_{p}=\alpha_{p}^{2}$. 
Since $E_{\lambda}$ has complex multiplication, looking at the sets of
primes $p$ for which the coefficients $\alpha_{p}, \beta_{p}$ and
$\gamma_{p}$ equal $\pm ip^{1/2}$ we easily see that the other two elliptic curves have
complex multiplication by the same quadratic field  and so up to a
twist the three weight two forms coincide. In particular $E_{\lambda}$
and $E_{\mu}$ are isogenous. 

\section{Involutions}
\label{sec:Inv}

In this section we will use an involution on a Calabi--Yau threefold
to split the cohomology group $H^{3}$. Note that van Geemen
and Nygaard (\cite{GeemenNygard}) were the first to use an automorphism
of a Calabi--Yau manifold to split the Galois representation and prove
modularity. 

\begin{prop}
Calabi--Yau threefolda $X_{53}$, $X_{244}$, $X_{267}$, $X_{274}$ and
$X_{275}$ are modular, with modular forms as listed in table \ref{tab1}.
\end{prop}

On some of the Calabi--Yau threefolds considered in this paper we can find
an involution. On the middle cohomology the involution may have only
eigenvalues $\pm1$. If both $1$ and $-1$ are eigenvalues then the
map gives us a splitting of $H^{3}$. Since the spliting is compatible
with the Frobenius morphism it is in fact a splitting of the Galois
representation into two--dimensional subrepresentations. 

We can use the Lefschetz formula to compute the trace of Frobenius
composed with the involution. This trace is equal to the trace of
Frobenius on the $+1$--eigenspace minus the trace of Frobenius on
the $-1$--eigenspace. Together with the trace of Frobenius on
$H^{3}$ this gives the traces on the two subspaces.

Assume that we have a $\mathbb Q$--linear involution on $\PP$ which
preserves the arrangement of eight planes.
This map induces an involution  $\Phi:X\lra X$ on the Calabi--Yau
threefold $X$ defined by this arrangement. We will compute
the trace
\[
d_{p}=\tr((\Frob_{p}\circ\Phi)^*|H^{3}(\bar X_{p}, \mathbb Q_{l}))
\]
of Frobenius composed with $\Phi$. Since this map acts by
multiplication with $\pm p$ on $H^{2}$ and with $\pm p^{2}$ on $H^{4}$ the
Lefschetz fixed--point formula relates $d_{p}$ to the number $N_{p}$
of fixed points of $\Frob_{p}\circ \Phi$.

\begin{lemma} 
  If $\Phi$ is a linear involution on $\PP[N](\bar{\mathbb F}_{p})$
  defined over $\mathbb F_{p}$ then
  the fixed points of $\Frob_{p}\circ\Phi$ are $\mathbb F_{p^{2}}$-rational. 
\end{lemma}
\begin{proof}
The Frobenius morphism $\Frob_{p}$ commutes with any linear involution
defined over $\mathbb F_{p}$, so any fixed point of $\Frob_{p}\circ\Phi$ is
also a fixed point of $\Frob_{p^{2}}$.
\end{proof}

Using the Lemma we reduce the counting of fixed points over
the infinite field $\bar{\mathbb F}_{p}$ to counting of points over the
finite field $\mathbb F_{p^{2}}$, which can easily be done using a computer.

From the representation as a Kummer fibration we can easily recognize
some linear involutions preserving the arrangement: 

\textbf{Arr. no. \phantom{0}53:} 
\((x,y,z,t)\mapsto(y,x,-t,-z)\) 

\textbf{Arr. no. 244:} 
\((x,y,z,t)\mapsto(y,x,-t,-z)\) 

\textbf{Arr. no. 267:} 
\((x,y,z,t)\mapsto(t,-z,-y,x)\) 

\textbf{Arr. no. 274:} 
\((x,y,z,t)\mapsto(z,-t,x,-y)\)

Simple computations show that the above involutions are not equal
to identity on the deformation space $H^{1}(\mathcal T_{X})\cong
H^{12}(X)$, hence they split the Galois
representations. In fact it is easy to observe that $H^{12}(X)\oplus
H^{21}(X)$ must be $(-1)$--eigenspaces. Counting fixed points on the
singular double octic yields for all primes $5\le p\le 97$:
\[\def\arraystretch{1.9}
\begin{array}[t]{r@{\hspace{3mm}}l}
X_{53}:&\hspace{3.4mm} 1+p^{3}-a_{p}+pb_{p}+p^{2}+p\\
X_{244}:&\def\arraystretch{1.1}
\begin{cases}
  1+p^{3}-a_{p}+pb_{p}+2p^{2}-p,&\qquad p\equiv 1\mod4\\
  1+p^{3}-a_{p}+pb_{p}+3p,&\qquad p\equiv 3\mod4
\end{cases}\\
X_{267}:&\hspace{3.4mm}1+p^{3}-a_{p}+pb_{p}+p^{2}-p\\
X_{274}:&
\begin{cases}
  1+p^{3}-a_{p}+pb_{p}+p^{2}-p,&\qquad p\equiv 1\mod4\\
  1+p^{3}-a_{p}+pb_{p}+p^{2}+3p,&\qquad p\equiv 3\mod4
\end{cases}
\end{array}
\]
Analyzing the action of Frobenius on the generators of the Picard
group and the space of curves $H^{4}$ gives the traces of Frobenius of
the two--dimensional Galois subrepresentations. Applying the
Faltings--Serre-Livn\'e method finishes the proof.

M. Sch\"utt suggested to us that counting points in $\mathbb F_{p}$ and
$\mathbb F_{p^{2}}$ we can compute the characteristic polynomial,
which factors into two degree two polynomials. Since we know that the
representation splits we get the traces of both actions. It is
however not straightforward that the numbers $a_{p}$ resp. $pb_p$
will correspond to the $+1$--eigenspace resp. the $-1$--eigenspace.

\begin{remark}
The described involutions act on singular double octics. Since the
resolution of singularities of a double octic is not unique (it depends
on the order in which we blow up lines in a triple point) it may
happen that an involution maps to a birational Calabi--Yau
threefold. Since two smooth models differ by a sequence of flops, we can
compose the involution with these flops or we can consider a threefold
that dominates both smooth models. The action on $H^{3}$ is well
defined. 

If we know that we can chose such a resolution of singularities of the
double covering to which the involution lifts, then the quotient will be
(after resolution) a rigid Calabi--Yau threefold.
\end{remark}

\begin{example}\label{ex:arr287}
Consider the arrangement of planes (arr. no. 287 in \cite{Meyer}) given by
\[
xyzt(x+y+z-3t)(x+y-3z+t)(x-3y+z+t)(-3x+y+z+t)=0.
\]

The corresponding Calabi--Yau threefold $X_{287}$ has Hodge numbers
$h^{11}(X_{287})=37$, $h^{12}(X_{287})=3$.
Counting points in $\mathbb F_{p}$ shows that, for $5\le p\le97$,
the trace of Frobenius on the middle cohomology equals $a_{p}+3b_{p}$,
where $a_{p}$ resp. $b_{p}$ are the coefficients of the weight 4 level 6
resp. weight 2 level 24 cusp form. The arrangement has many linear
symmetries. We can use the induced involutions on $X$ to decompose
the Galois representation.

We can also use the elliptic fibrations on $X$ described in
\cite{Meyer} and apply the deformation argument from example~\ref{ex:ar269}
to prove modularity of $X_{287}$.

In fact the full permutation group $S_{4}$ acts on this Calabi--Yau
threefold. If we consider the action of permutations of order 3, then
the eigenvalues will be defined in $\mathbb F_{p}$ only for some $p$,
so the decomposition of Frobenius action will depend on $p$.
\end{example}

\subsection*{Acknowledgements}
The work on this paper was done during the first named author's stays at the
Institutes of Mathematics of the Johannes Gutenberg-Universit\"at Mainz
and the Universit\"at Hannover. He would like to thank both institutions
for their hospitality. The authors also would like to thank Prof. Duco van Straten,
Prof. Klaus Hulek and Matthias Sch\"utt for their help.

\end{document}